# Structure and characterization of ruled surfaces in Minkowski 3-space

Fatma Güler, Emin Kasap
Department of Mathematics, Arts and Science Faculty,
Ondokuz Mayis University, 55139 Samsun, Turkey.
e-mail: f.guler@omu.edu.tr
e-mail: kasape@omu.edu.tr

**Abstract**

In this paper, we consider non developable ruled surface with spacelike ruling, timelike ruling, respectively. We give the relations between the structure functions with the curvature and torsion of the striction line of the timelike and spacelike non developable ruled surfaces. Also, we have calculated the gaussian and mean curvatures of timelike and spacelike non developable ruled surfaces using the structure functions.

**Introduction**

It is known that a closed ruled surface has two integral invariants in called the pitch and the angle of pitch,[4], [5]. In recent years several authors have used these invariants in their investigations in [1], [2]. In one of the latest papers, in [3] the authors studied the relations between the structure functions of the ruled surface and the curvature, torsion of the striction line of the ruled surface. Also, their introduced the invariants of non developable ruled surfaces and kinematical characterizations of non developable ruled surface in Euclidean 3-space. In [2] , the authors defined pitch function for any non developable ruled surfaces and use this notion to give a new characterization of $B$ scrolls in Minkowski 3- space. In this paper, we give the relations between the structure functions of the timelike and spacelike non developable surfaces and the curvature, torsion of the striction line of the timelike and spacelike non developable ruled surfaces in Minkowski 3-space. Then, we study timelike and spacelike non developable ruled surfaces in Minkowski 3-space in terms of the mean curvature and the gaussian curvature.

**Preliminaries**

Let us consider Minkowski 3-space $IR_1^3 = \left[ IR_1^3, (+,+,-) \right]$ and let the Lorentzian iner product of $X = (x_1, x_2, x_3)$ and $Y = (y_1, y_2, y_3) \in IR_1^3$ be $\langle X, Y \rangle = x_1 y_1 + x_2 y_2 - x_3 y_3$. A vector $X \in IR_1^3$ is called a spacelike vector when $\langle X, X \rangle > 0$ or $X = 0$. It is called timelike and null(ligtlike) vector in case of $\langle X, X \rangle < 0$ and $\langle X, X \rangle = 0$ for $X \neq 0$, respectively,[6].

The vector product of vectors $X = (x_1, x_2, x_3)$ and $Y = (y_1, y_2, y_3)$ in $IR_1^3$ is defined by [7],

$$X \times Y = (x_2 y_3 - x_3 y_2, x_1 y_3 - x_3 y_1, x_2 y_1 - x_1 y_2).$$

A regular curve $\alpha(s): I \to IR_1^3$, $I \in IR$ in $IR_1^3$ is said to be a spacelike, timelike and null curve if the velocity vector $\alpha'(s) = d\alpha/ds$ is a spacelike, timelike or null vector, respectively,[8].

A surface in $IR_1^3$ is called a timelike surface if the normal vector on the surface spacelike vector and is called spacelike surface if the normal vector on the surface timelike vector.

Let $X(u,v) = \alpha(u) + vb(u)$ be a non developable ruled surface in $E^3$ with $b^2(u) = 1$ and the parameter $u$ is the arc length parameter of $(b(u))$ as a unit spherical curve in $E^3$. The spherical Frenet formulas of $(b(u))$ unit spherical curve can be written as

$$x'(u) = \alpha(u)$$
$$\alpha'(u) = -x(u) + k_g(u) y(u)$$
$$y'(u) = -k_g(u) \alpha(u)$$

where $x(u) = b(u)$, $x'(u) = \alpha(u)$ and $y(u) = \alpha(u) \times x(u)$, [3].

The pitch $\delta(u_0)$ of the ruled surface $X(u,v) = \alpha(u) + vb(u)$ with spacelike ruling or timelike ruling at $\alpha(u_0)$ is defined by

$$\delta(u_0) := \lim_{\Delta u \to 0} \frac{[A(u_0 + \Delta u) - \alpha(u_0 + \Delta u)].b(u_0 + \Delta u)}{\Delta u} = -\langle \alpha'(u_0), b(u_0) \rangle$$

where $\delta(u)$ the pitch function of the ruled surface $X(u,v)$ in $IR_1^3$, [2].

For one parameter unit vector field $b(u)$ with $|b'(u)| = 1$, the function

$$\theta(u) := -\langle (b'(u) \times b(u))', b'(u) \rangle$$

is called angle (density) function or called self spinning (density) function of vector field $b(u)$, [1].

Let $X(u,v) = \alpha(u) + vb(u)$ be any non developable ruled surface in $E^3$ and $\alpha(u)$ the striction line of $X(u,v)$ such that $|b(u)| = |b'(u)| = 1$. The angle function (self spinning function) of vector field $b(u)$ is called angle (density) function of pitch or self spinning (density) function of the non developable ruled surface $X(u,v)$, [1]. Therefore, we have $\theta(u) = k_g(u)$.

Then the angle angle (density) function of pitch or self spinning (density) function of one parameter unit vector field $b(u)$ is the spherical curvature function of which $b(u)$ is considered as a unit spherical curve in $E^3$, [1].

Let $X(u,v) = \alpha(u) + vb(u)$ be any non developable ruled surface and $\alpha(u)$ the striction line of $X(u,v)$ such that $\alpha'(u) = \lambda(u) x(u) + \mu(u) y(u)$, here $\{\alpha(u), x(u) = b(u), y(u)\}$ is the spherical Frenet frame of the spherical curve $(b(u))$, [3].

## 3. Non Developable Timelike Ruled Surface in Minkowski 3-Space

Let $X(u,v) = \alpha(u) + vb(u)$ be a non developable timelike ruled surface in Minkowski 3-space $IR_1^3$ with $\langle b(u), b(u) \rangle = \mp 1$ and $|\langle b'(u), b'(u) \rangle| = 1$, that is, $u$ is the arc length parameter of $b(u)$ as a curve on the Sitter space $S_1^2$ (when $\langle b(u), b(u) \rangle = 1$) or hyperbolic space $H^2$ (when $\langle b(u), b(u) \rangle = -1$). We also assume that the base curve $\alpha(u)$ of the timelike ruled surface $X(u,v)$ is the striction line of the surface, that is $\langle \alpha'(u), b'(u) \rangle = 0$.

### 3.1. Non Developable Timelike Ruled Surface with Spacelike Direction

Let $X(u,v) = \alpha(u) + vb(u)$ be a non developable timelike ruled surface with $\langle b(u), b(u) \rangle = 1$. Choosing $x(u) = b(u)$, $x'(u) = \alpha(u)$ and $y(u) = \alpha(u) \times x(u)$, the spherical Frenet formulas of unit spherical curve $(b(u))$ can be written as,

$$x'(u) = \alpha(u)$$
$$\alpha'(u) = \varepsilon x(u) - k_g(u) y(u)$$
$$y'(u) = -k_g(u) \alpha(u)$$

where $\varepsilon = \begin{cases} 1, & \alpha(u) \text{ timelike vector} \\ -1, & \alpha(u) \text{ spacelike vector.} \end{cases}$

We consider the properties and relations of the structure functions of non developable timelike ruled surface with spacelike direction.

Let $X(u,v) = \alpha(u) + vb(u)$ be a non developable timelike ruled surface with spacelike direction(when $\langle \alpha(u), \alpha(u) \rangle = -1$) and $\alpha(u)$ the striction line of $X(u,v)$ such that $\alpha'(u) = -\lambda(u) x(u) + \mu(u) y(u)$, here $\{\alpha(u), x(u), y(u)\}$ is the spherical frenet frame of the spherical curve $(b(u))$, $u$ is the arc length parameter of $b(u)$. The definition of the structure function $\lambda(u)$. We have $\lambda(u) = -\delta(u)$.

**Theorem 3.1.1.**

Then the curvature function $\kappa(u)$ and torsion function $\tau(u)$ of the striction line $\alpha(u)$ of $X(u,v)$ are given by

$$\kappa^2 = \frac{(\lambda+\mu k_g)^2(\lambda^2+\mu^2)+(\mu\lambda'-\lambda\mu')^2}{(\lambda^2+\mu^2)^3}$$

$$\tau = \frac{(\lambda+\mu k_g)(\lambda^2 k_g - \lambda\mu k_g^2 + \lambda\mu'' - \lambda'' + \lambda + \mu k_g) + (2\lambda' + k_g'\mu + 2k_g\mu')(\mu\lambda' - \lambda\mu')}{(\lambda^2+\mu^2)(\lambda+\mu k_g)^2 + (\mu\lambda'-\lambda\mu')^2}$$

Also, the Gauss curvature and mean curvature of $X(u,v)$ are

$$K(u,v) = -\frac{\mu^2}{\mu^2-v^2}, \quad H(u,v) = \frac{-\mu(\lambda+\mu k_g) - v(\mu'-vk_g+2\lambda\mu)}{2(\mu^2-v^2)}.$$

**Proof.** From

$$\kappa(u) = \frac{\|\alpha'(u)\times\alpha''(u)\|}{\|\alpha'(u)\|^3}, \quad \tau(u) = \frac{\langle\alpha'(u),\alpha''(u)\times\alpha'''(u)\rangle}{\|\alpha'(u)\times\alpha''(u)\|^2}$$

and

$$K(u,v) = -\varepsilon\frac{LN-M^2}{EG-F^2}, \quad H(u,v) = \varepsilon\frac{GL+EN-2FM}{2(EG-F^2)}.$$

By a direct calculation we can get the conclusion of this theorem, where

$$\varepsilon = \begin{cases} -1, & n \text{ timelike vector} \\ +1, & n \text{ spacelike vector} \end{cases}.$$

The unit normal vector $n$ spacelike vector

$$n = \frac{\mu(u)\alpha(u)+vy(u)}{\sqrt{\mu^2(u)+v^2}}$$

where $\lambda(u)$, $\mu(u)$ and $k_g(u)$ structure functions of $X(u,v)$. Also, the first fundamental quantities of $X(u,v)$ are

$$E = \lambda^2(u) + \mu^2(u) - v^2, \quad F = -\lambda(u), \quad G = 1.$$

The second fundamental quantities are

$$L = -\mu(\lambda+\mu k_g) - v(\mu'-vk_g), \quad M = \mu, \quad N = 0.$$

The Gauss curvature and mean curvature of $X(u,v)$ are

$$K(u,v) = -\frac{\mu^2}{\mu^2-v^2}, \quad H(u,v) = \frac{-\mu(\lambda+\mu k_g) - v(\mu'-vk_g+2\lambda\mu)}{2(\mu^2-v^2)}.$$

Let $X(u,v) = \alpha(u) + vb(u)$ be a non developable timelike ruled surface with spacelike direction(when $\langle\alpha(u),\alpha(u)\rangle = +1$) and $\alpha(u)$ the striction line of $X(u,v)$ such that $\alpha'(u) = \lambda(u)x(u) + \mu(u)y(u)$.

**Theorem 3.1.2.**

Then the curvature function $\kappa(u)$ and torsion function $\tau(u)$ of the striction line $\alpha(u)$ of $X(u,v)$ and the Gauss curvature $K(u,v)$ and mean curvature of $H(u,v)$ of $X(u,v)$ are given by

$$\kappa^2 = \frac{(\mu k_g - \lambda)^2 (\lambda^2 + \mu^2) + (\lambda\mu' - \mu\lambda')^2}{(\lambda^2 + \mu^2)^3}$$

$$\tau = \frac{(\lambda - \mu k_g)(-\lambda^2 k_g + \lambda\mu k_g^2 - \lambda\mu'' - \lambda'' + \lambda - \mu k_g) + (2\lambda' - k_g'\mu - 2k_g\mu')(\mu\lambda' - \lambda\mu')}{(\lambda^2 + \mu^2)(\mu k_g - \lambda)^2 + (\lambda\mu' - \mu\lambda')^2}$$

and

$$K(u,v) = -\frac{\mu^2}{\mu^2 + v^2 - 2\lambda^2}, \quad H(u,v) = \frac{\lambda v - \mu\mu' + 2\lambda\mu}{2(\mu^2 + v^2 - 2\lambda^2)}.$$

**Proof.** The first fundamental quantities of $X(u,v)$ are

$$E = \mu^2(u) - \lambda^2(u) + v^2, \quad F = -\lambda(u), \quad G = 1.$$

The unit normal vector is spacelike vector

$$n = \frac{-\mu(u)\alpha(u) + vy(u)}{\sqrt{\mu^2(u) + v^2}}.$$

The second fundamental quantities are

$$L = \lambda(u)v - \mu(u)\mu'(u), \quad M = \mu(u), \quad N = 0.$$

The Gauss curvature and mean curvature of $X(u,v)$ are

$$K(u,v) = -\frac{\mu^2}{\mu^2 + v^2 - 2\lambda^2}, \quad H(u,v) = \frac{\lambda v - \mu\mu' + 2\lambda\mu}{2(\mu^2 + v^2 - 2\lambda^2)}.$$

### 3.2. Non Developable Timelike Ruled Surface with Timelike Direction

Let $X(u,v) = \alpha(u) + vb(u)$ be a non developable timelike ruled surface with $\langle b(u), b(u) \rangle = -1$. The spherical Frenet formulas of unit spherical curve $(b(u))$ can be written as,

$$x'(u) = \alpha(u)$$
$$\alpha'(u) = -x(u) - k_g(u) y(u)$$
$$y'(u) = -k_g(u)\alpha(u)$$

where $x(u) = b(u)$, $x'(u) = \alpha(u)$ and $y(u) = \alpha(u) \times x(u)$.

Let $X(u,v) = \alpha(u) + vb(u)$ be a non developable timelike ruled surface with timelike direction and $\alpha(u)$ the striction line of $X(u,v)$ such that $\alpha'(u) = -\lambda(u)x(u) + \mu(u)y(u)$. The definition of the structure function $\lambda(u)$. We have $\lambda(u) = -\delta(u)$.

**Theorem 3.2.2.**

Then the curvature function $\kappa(u)$ and torsion function $\tau(u)$ of the striction line $\alpha(u)$ of $X(u,v)$ are given by

$$\kappa^2 = \frac{(\lambda + \mu k_g)^2 (\lambda^2 + \mu^2) + (\mu\lambda' - \lambda\mu')^2}{(\lambda^2 + \mu^2)^3}$$

$$\tau = \frac{(\lambda + \mu k_g)(\lambda^2 k_g + \lambda \mu k_g^2 + \lambda\mu'' - \lambda'' + \lambda + \mu k_g) + (2\lambda' - k_g'\mu + 2k_g\mu')(\mu\lambda' - \lambda\mu')}{(\lambda^2 + \mu^2)(\lambda + \mu k_g)^2 + (\mu\lambda' - \lambda\mu')^2}$$

Also, the Gauss curvature and mean curvature of $X(u,v)$ are

$$K(u,v) = \frac{\mu^2(u)}{\mu^2(u) + v^2}, \quad H(u,v) = \frac{\mu(\lambda - \mu k_g) - v(\mu' - vk_g - 2\lambda\mu)}{2(\mu^2 + v^2)}.$$

**Proof.**

The first fundamental quantities of $X(u,v)$ are
$$E = \mu^2(u) - \lambda^2(u) + v^2, \quad F = -\lambda(u), \quad G = -1.$$
The unit normal vector is timelike vector
$$n = \frac{-\mu(u)\alpha(u) + vy(u)}{\sqrt{\mu^2(u) + v^2}}.$$
The second fundamental quantities are
$$L = \mu(u)(\lambda(u) - \mu(u)k_g(u)) - v(\mu'(u) - vk_g(u)), \quad M = \mu(u), \quad N = 0.$$

The Gauss curvature and mean curvature of $X(u,v)$ are

$$K(u,v) = \frac{\mu^2(u)}{\mu^2(u) + v^2}, \quad H(u,v) = \frac{\mu(\lambda - \mu k_g) - v(\mu' - vk_g - 2\lambda\mu)}{2(\mu^2 + v^2)}.$$

## 4. Non Developable Spacelike Ruled Surface in Minkowski Space

Let $X(u,v) = \alpha(u) + vb(u)$ be a non developable spacelike ruled surface. Choosing $x(u) = b(u)$, $x'(u) = \alpha(u)$ and $y(u) = \alpha(u) \times x(u)$, the spherical Frenet formulas of unit spherical curve $(b(u))$ can be written as,

$$x'(u) = \alpha(u)$$
$$\alpha'(u) = \varepsilon x(u) - k_g(u) y(u)$$
$$y'(u) = -k_g(u) \alpha(u)$$

where $\varepsilon = \begin{cases} 1, & y(u) \text{ spacelike vector} \\ -1, & y(u) \text{ timelike vector.} \end{cases}$

**Theorem 4.1.**

Then the curvature function $\kappa(u)$ and torsion function $\tau(u)$ of the striction line $\alpha(u)$ of $X(u,v)$ and the Gauss curvature $K(u,v)$ and mean curvature of $H(u,v)$ of $X(u,v)$ are given by

$$\kappa^2 = \frac{(\lambda + \mu k_g)^2 (\lambda^2 + \mu^2) + (\mu\lambda' - \lambda\mu')^2}{(\lambda^2 + \mu^2)^3}$$

$$\tau = \frac{(\lambda + \mu k_g)(\lambda^2 k_g - \lambda\mu k_g^2 + \lambda\mu'' - \lambda'' + \lambda + \mu k_g) + (2\lambda' - k_g'\mu + 2k_g\mu')(\mu\lambda' - \lambda\mu')}{(\lambda^2 + \mu^2)(\lambda + \mu k_g)^2 + (\mu\lambda' - \lambda\mu')^2}$$

and

$$K(u,v) = \frac{\mu^2(u)}{\mu^2(u) - v^2}, \quad H(u,v) = \frac{\mu(\lambda + \mu k_g) + v(\mu' - vk_g - 2\lambda\mu)}{2(\mu^2 - v^2)}$$

where $y(u)$ timelike vector,

$$K(u,v) = \frac{\mu^2(u)}{\mu^2(u) + v^2 - 2\lambda^2(u)}, \quad H(u,v) = \frac{\mu(u)\mu'(u) - \lambda(u)v - 2\lambda(u)\mu(u)}{2(\mu^2(u) + v^2 - 2\lambda^2(u))}$$

where $y(u)$ spacelike vector.

**Proof.**

The first fundamental quantities of $X(u,v)$ are
$$E = \mu^2(u) + \lambda^2(u) - v^2, \quad F = -\lambda(u), \quad G = 1.$$

The unit normal vector is timelike vector
$$n = \frac{\mu(u)\alpha(u) + vy(u)}{\sqrt{\mu^2(u) + v^2}}.$$

The second fundamental quantities are
$$L = -\mu(u)\big(\lambda(u) + \mu(u)k_g(u)\big) - v\big(\mu'(u) - vk_g(u)\big), \quad M = \mu(u), \quad N = 0.$$

The Gauss curvature and mean curvature of $X(u,v)$ are
$$K(u,v) = \frac{\mu^2(u)}{\mu^2(u) - v^2}, \quad H(u,v) = \frac{\mu(\lambda + \mu k_g) + v(\mu' - vk_g - 2\lambda\mu)}{2(\mu^2 - v^2)}$$

where $y(u)$ timelike vector.

The first fundamental quantities of $X(u,v)$ are
$$E = \mu^2(u) - \lambda^2(u) + v^2, \quad F = -\lambda(u), \quad G = 1.$$

The unit normal vector is timelike vector
$$n = -\frac{\mu(u)\alpha(u) + vy(u)}{\sqrt{\mu^2(u) + v^2}}.$$

The second fundamental quantities are
$$L = \lambda(u)v - \mu(u)\mu'(u), \quad M = \mu(u), \quad N = 0.$$

The Gauss curvature and mean curvature of $X(u,v)$ are
$$K(u,v) = \frac{\mu^2(u)}{\mu^2(u) + v^2 - 2\lambda^2}, \quad H(u,v) = \frac{\mu(u)\mu'(u) - \lambda(u)v - 2\lambda(u)\mu(u)}{2(\mu^2(u) + v^2 - 2\lambda^2(u))}$$

where $y(u)$ spacelike vector.

# References

[1] Huili Liu, Yanhua Yu, Seoung Dal Jung, Invariants of non-developable ruled surfaces in Euclidean 3-space, Preprint, 2011, NEU.

[2] Huili Liu, Yuan Yuan, Pitch functions of ruled surfaces and B-scrolls in Minkowski 3-space, Journal of Geometry and Physics, 62(2012), 47-52.

[3] Yanhua Yu, Huili Liu, Seoung Dal Jung, Structure and Characterization of Ruled Surfaces in Euclidean 3-Space, Applied Mathematics and Computation,2014.

[4] H.R.Müller, Monatsh Math, 53(1951), pp.206-214.

[5] J. Hoschek Avch. Math., XXIV(1973), pp.218-224.

[6] B. O'Neill, Semi- Riemannian Geometry, Academic Pres, New York, 1983.

[7] K. Akutagawa, S. Nishikawa, The Gauss map and spacelike surfaces with prescribed mean curvature in Minkowski 3-space, Tohoku Math. J. 42(1990)67-82.

[8] T.Ikaw, On curves and sub manifolds in indefinite Riemannian manifold, Tsukuba J.Math. 9(2)(1985) 353-371.